\documentclass[12pt]{article}
\usepackage{a4wide}
\usepackage{epsfig}
\usepackage{amsfonts}
\usepackage{amsthm}
\usepackage{amsmath}
\usepackage{amssymb}
\newtheorem{theorem}{Theorem}
\newtheorem{lemma}[theorem]{Lemma}
\def\const{{0.299871}}
\def\constsup{{0.299309}}
\def\constfive{{0.313972}}
\def\PP{{\cal P}}
\def\RR{{\cal R}}
\providecommand{\norm}[1]{\lVert#1\rVert}

\def\bx#1{\hbox to 10pt{\hfil$#1$\hfil}}
\def\.{\hbox{$\bullet$}}
\def\+{\hbox{\tt+}}
\def\x{\hbox{$\times$}}
\def\o{\hbox{$\circ$}}
\def\?{\hbox{\tt?}}

\begin{document}
\title{Domination number of cubic graphs with large girth\thanks{This research was supported by the Czech-Slovenian bilateral project MEB 090805 (on the Czech side) and BI-CZ/08-09-005 (on the Slovenian side).}}
\author{Daniel Kr\'al'\thanks{Institute for Theoretical Computer Science (ITI), Faculty of Mathematics and Physics, Charles University, Malostransk\'e n\'am\v est\'\i{} 25, 118 00 Prague 1, Czech Republic. E-mail: {\tt kral@kam.mff.cuni.cz}. Institute for Theoretical computer science is supported as project 1M0545 by Czech Ministry of Education. This research was also supported by the grant GACR 201/09/0197.}\and
        Petr {\v S}koda\thanks{Department of Applied Mathematics, Faculty of Mathematics and Physics, Charles University, Malostransk\'e n\'am\v est\'\i{} 25, 118 00 Prague 1, Czech Republic. E-mail: {\tt peskoj@kam.mff.cuni.cz}.}\and
	Jan Volec\thanks{Department of Applied Mathematics, Faculty of Mathematics and Physics, Charles University, Malostransk\'e n\'am\v est\'\i{} 25, 118 00 Prague 1, Czech Republic. E-mail: {\tt janv@kam.mff.cuni.cz}.}}
\date{}
\maketitle
\begin{abstract}
We show that every $n$-vertex cubic graph with girth at least $g$ have domination
number at most $\const{}n+O\left(n/g\right)<3n/10+O\left(n/g\right)$.
\end{abstract}

\section{Introduction}

The notion of a dominating set is a classical notion in graph theory with a large amount of literature associated with it.
For the sake of completeness, let us recall that a set $D$ of vertices of a graph $G$ is {\em dominating}
if every vertex of $G$ is contained in $D$ or has a neighbor in $D$ and the {\em domination number $\gamma(G)$} of $G$
is the smallest size of a dominating set of $G$. In this paper, we study the domination number of cubic graphs,
i.e., graphs where every vertex has degree three.

In 1996, Reed~\cite{bib-reed} conjectured that every $n$-vertex connected cubic graph has domination number
at most $\left\lceil n/3\right\rceil$. Though the conjecture turn out to be false~\cite{bib-kostochka+,bib-kelmans},
the conjecture becomes true with additional assumption that the cubic graph has {\em girth} at least $g$,
i.e., it has no cycles of length less than $g$.
The first results in this direction are the bounds on the domination number of $\left(\frac{1}{3}+\frac{1}{3g+3}\right)n$
of Kawarabayashi et al.~\cite{bib-kawarabayashi+} for bridgeless $n$-vertex cubic graphs
with girth at least $g$ for $g$ divisible by three and
$\left(\frac{1}{3}+\frac{8}{3g^2}\right)n$ of Kostochka and Stodolsky~\cite{bib-kostochka+bound}
for all $n$-vertex cubic graphs with girth at least $g$. The magic threshold of $n/3$ was first beaten
for cubic graph with large girth by L{\"o}wenstein and Rautenbach~\cite{bib-lowenstein+}
who showed that every $n$-vertex cubic graph with girth at least $g\ge 5$
contains a dominating set of size at most $\left(\frac{44}{135}+\frac{82}{135g}\right)n\approx 0.325926n+O\left(n/g\right)$
The bound was further improved by Rautenbach and Reed~\cite{bib-rautenbach+} to $0.321216n+O\left(n/g\right)$.

We further improve these bounds and manage to lower them below the $3n/10$ threshold. Our main result
is the following:

\begin{theorem}
\label{thm-main}
Let $G$ be an $n$-vertex cubic graph with girth at least $g$. The domination number of $G$
is at most $\const{}n+O\left(n/g\right)\le 3n/10+O\left(n/g\right)$.
\end{theorem}

\noindent At this point, we remark that numerical computations involved in our argument were done using
a computer (but the rules given in Figure~\ref{fig-rules} were generated by hand)
though the whole proof can be easily verified to be correct without computer assistance.

Before we start the exposition of our proof, let us mention a connection to random cubic graphs.
A random cubic graph almost surely contains only a bounded number of cycles of length less than $g$
for every fixed integer $g$ and thus the value of the domination number of a random cubic graph
should indicate how tight our result could be. It is known that the domination number $\gamma$ of 
a random cubic graph is almost surely at least $0.2636n$~\cite{bib-molloy+} and
at most $0.2794n$~\cite{bib-duckworth+} which should indicate space for further improvements.

\section{Proof}

We first provide a general overview of our method and illustrate it on a small example.
At the end, we then apply the method in a setting yielding Theorem~\ref{thm-main}.

\subsection{Overview}
\label{sub-overview}

In this subsection,
we explain main ideas of our method and we later provide necessary technical details related to it.
Fix a cubic bridgeless graph $G$ with girth at least $g$ and also fix an integer $K$
which determines the number of levels as defined later.
Consider a $2$-factor of $G$ (which exists by the Petersen theorem) and
decompose each cycle of the $2$-factor into vertex-disjoint paths
$P_1,\ldots,P_{\ell}$ with the number of vertices between $g/4K$ and $g/2K$ (this is possible
since the length of each cycle of the $2$-factor is at least $g$).
The vertices of the paths $P_1,\ldots,P_{\ell}$ are considered to be ordered from one end of the path
towards the other; a {\em mate} of a vertex $v$ is the neighbor of $G$ not adjacent to $v$
on the cycle of the $2$-factor.

A dominating set $D$ of $G$ will be given by a labeling of vertices of $G$ we construct.
Each vertex of $G$ will be assigned an input label
which is one of the symbols $\+$, $\x$, $\.$ and $\o$, and
an output label which is one of the symbols $\oplus$, $\otimes$ and $\odot$.
The dominating set $D$ will contain the vertices with input label $\+$ or output label $\otimes$ (as well as
several others, see the next subsection for an exhaustive definition).
We explain the intuitive meaning of the labels later.

Split the paths $P_1,\ldots,P_k$ into $K$ sets $\PP_1,\ldots,\PP_K$
including each path to a single set randomly uniformly and independently of the other paths.
The sets $\PP_1,\ldots,\PP_K$ are referred to as levels and vertices on paths in $\PP_i$
are said to be on the level $i$.
First, the vertices contained in paths of $\PP_1$ are assigned labels, then those in paths of $\PP_2$, etc.
Let $P$ be a path included in $\PP_i$ and assume that the vertices of the paths in $\PP_1\cup\cdots\cup\PP_{i-1}$
have already been labelled. The input label of a vertex of $P$ is:
\begin{itemize}
\item the symbol $\+$ if its mate is on a path in $\PP_1\cup\cdots\cup\PP_{i-1}$ and its output label is $\oplus$,
\item the symbol $\x$ if its mate is on a path in $\PP_1\cup\cdots\cup\PP_{i-1}$ and its output label is $\otimes$,
\item the symbol $\.$ if its mate is on a path in $\PP_1\cup\cdots\cup\PP_{i-1}$ and its output label is $\odot$,
\item the symbol $\.$ if its mate is on a path in $\PP_i$, and
\item the symbol $\o$ if its mate is on a path in $\PP_{i+1}\cup\cdots\cup\PP_k$.
\end{itemize}
The output labels are assigned in blocks using rules.
Each rule is a pair of a sequence of input symbols and the symbol $\?$ which represents a wild-card and
a sequence of output symbols. The lengths of the two sequences will always be the same.
We usually write an arrow between the two sequences, e.g.,
one of the rules can be $\+\?\rightarrow\odot\odot$. Naturally, the sequence of input symbols and $\?$
is called the left-hand side of the rule and the sequence of output symbols the right-hand side.
Finally, if $\sigma\rightarrow\tau$ is a rule, we use
$\sigma_i$ is the $i$-th symbol of $\sigma$ and $\tau_i$ is the $i$-th symbol of $\tau$.
An example of a set of rules is given in Figure~\ref{fig-example}.

\begin{figure}
\begin{center}
$$\begin{array}{lclclcl}
\bx\+\bx\?  & \rightarrow &  \bx\odot\bx\odot &\quad& \bx\.\bx\o\bx\?  & \rightarrow &  \bx\odot\bx\otimes\bx\odot\\
\bx\x  & \rightarrow &  \bx\odot && \bx\o\bx\+\bx\?  & \rightarrow &  \bx\odot\bx\odot\bx\odot\\
\bx\.\bx\+\bx\?  & \rightarrow &  \bx\odot\bx\odot\bx\odot && \bx\o\bx\x  & \rightarrow &  \bx\oplus\bx\odot\\
\bx\.\bx\x\bx\?  & \rightarrow &  \bx\odot\bx\otimes\bx\odot && \bx\o\bx\.\bx\?  & \rightarrow &  \bx\odot\bx\otimes\bx\odot\\
\bx\.\bx\.\bx\?  & \rightarrow &  \bx\odot\bx\otimes\bx\odot && \bx\o\bx\o\bx\?  & \rightarrow &  \bx\odot\bx\otimes\bx\odot%
\end{array}$$
\end{center}
\caption{An example of a (correct) set of rules.}
\label{fig-example}
\end{figure}

We look for a rule whose left-hand side matches the input labels of vertices at the beginning of $P$.
All the considered rule sets will have the property that such a rule is unique.
The right-hand side of the rule then determine the output labels. We then move after the vertices
with output labels assigned and again look for a rule whose left-hand side matches the input labels
of vertices without output labels. After finding a suitable rule, the next group of vertices
is assigned output labels, we move after them and continue in this way until we reach the end of the path $P$.
Note that it can happen that we are left with few vertices without output labels at the end of the path $P$---this
will be handled in the next subsection.

The meaning of output labels is the following: if a vertex $v$ is labelled with $\otimes$,
then $v$ should be included to $D$. If a vertex $v$ is labelled with $\oplus$, then its mate
should be included to $D$ (this output label can only be assigned to vertices with input labels $\o$,
i.e., with mates on higher levels).
Finally, if $v$ is labelled with $\odot$, then it is dominated by one of its neighbors on its path or
by its mate on a lower level.
The input labels can be interpreted as follows: the label $\+$ represents that $v$ is included
to $D$, the label $\x$ represents that its mate on a lower label is included,
the label $\.$ represents that neither $v$
nor its mate has yet not been included to $D$ and the mate is on a lower or the same level and
the label $\o$ represents that a mate is on a higher level.

A set of rules is called {\em correct} if for every rule $\sigma\rightarrow\tau$ the following holds:
\begin{itemize}
\item if $\sigma_i$ is $\.$ or $\o$ and neither $\sigma_{i-1}$ nor $\sigma_{i+1}$ (if they exist) is $+$,
      then $\tau_i$ is $\oplus$ or
      one of the symbols $\tau_{i-1}$, $\tau_i$ and $\tau_{i+1}$ is $\otimes$, and
\item if $\tau_i$ is $\oplus$, then $\sigma_i$ is $\o$.
\end{itemize}
Intuitively, if a set of rules is correct,
then the set $D$ containing the vertices with input label $\+$ or output label $\otimes$ is always dominating and
the label $\oplus$ can be only assigned to vertices with mates on higher levels.
The set of rules given in Figure~\ref{fig-example} is correct.

\subsection{Analysis and adjustments}
\label{sub-analysis}

In this subsection, we provide further details on the labeling procedure and analyze it.
The first thing to cope with is the fact that a cubic graph need not to have a $2$-factor.
This is handled in a way analogous to that used in~\cite{bib-rautenbach+}.
For a collection of vertex-disjoint paths, a vertex is {\em covered} if it is contained in one of the paths.

\begin{lemma}
\label{lm-path}
For every $K\ge 2$, $g$ and $n$-vertex cubic graph $G$ with girth at least $g$,
there exists a collection of at most $\frac{3+8K}{2g}n$ vertex disjoint paths in $G$
with the number of vertices less than $g/2K$ that covers at least $n-O\left(n/g\right)$ vertices of $G$.
Moreover, all vertices on the paths can be grouped into pairs of vertices
adjacent through an edge not contained in the paths.
\end{lemma}

\begin{proof}
We only sketch the proof as it essentially follows the lines of the argument in~\cite{bib-rautenbach+}.
Let $S$ be a set of vertices of $G$. The number of odd components of $G\setminus S$ with a single edge to $S$
is at most $n/g$ (since each such component contains a cycle; in fact their number can be bounded by $n/2^{\Theta(g)}$,
but we do not need this finer estimate here). Hence, the Tutte-Berge formula implies that $G$
has a matching $M$ covering all but at most $n/g$ vertices of $G$. Remove the vertices not covered by $M$ and
the edges of $M$. In this way, we obtain a collection of at most $3n/2g$ paths and cycles,
each of them with length at least $g$.

We split the cycles into paths with the number of vertices at least $g/4K$ and less than $g/2K$.
This is clearly possible. Similarly, the paths with at least $g/2K$ vertices are split.
Paths with less than $g/2K$ are preserved and there is at most $3n/2g$ such paths.
Altogether, the final collection of paths consists of at most $4Kn/g$ new paths (as each of them
has at least $g/4K$ vertices) and at most $3n/2g$ paths of the original collection.
\end{proof}

Fix a cubic graph $G$ and a collection of vertex disjoint paths $P_1,\ldots,P_k$ as described in Lemma~\ref{lm-path}.
We say that a vertex $w$ is {\em $1$-close} to another vertex $v$ if $w$ lies on the same path as
the mate of $v$ (in particular, $w$ can be the mate of $v$). Note that the relation of being ``$1$-close''
is not symmetric. For $i>1$, we say that $w$ is $i$-close to $v$ if $w$ is $1$-close to $v$ or
$w$ is $(i-1)$-close to a vertex $w'$ where $w'$ is not the mate of $v$ but lies on the same path as the mate of $v$.
Observe that if $w$ is $i$-close to $v$, then the distance between $v$ and $w$ is less than $gi/2K$.
Because of this, there is no vertex $w$ that is $K$-close to two different
vertices $v$ and $v'$ lying on the same path (if such $w$ existed, then the
paths from $v$ and $v'$ to $w$ witnessing that $w$ is $K$-close and the subpath
between $v$ and $v'$ would form a cycle of length less than $g$).
The notion of being $K$-close will play a crucial role in our argument that
certain events are independent later in the proof.

For each path $P$ among $P_1,\ldots,P_k$, choose randomly uniformly and
independently of other paths an integer $i$, $1\le i\le K$, and include $P$ in the set $\PP_i$.
We now recursively define the probabilities $p_{i}(\+)$, $p_{i}(\x)$, $p_{i}(\.)$ and $p_{i}(\o)$
which represent probabilities that a vertex of a path in $\PP_i$ has a certain input label,
the probabilities $q_{i}(\oplus)$, $q_{i}(\otimes)$ and $q_{i}(\odot)$ that it has a certain output label and
$q^{\o}_{i}(\oplus)$, $q^{\o}_{i}(\otimes)$ and $q^{\o}_{i}(\odot)$ that it has a certain output label conditioned
on input label being equal to $\o$.
These numbers represent the probabilities in the ``ideal'' case, which includes the assumption that
all the labels are ``independent'' of each other and all the paths in $\PP_i$ are infinite, and
we will have to argue that they can also be applied to our labeling procedure.

Let us start with estimating the input probabilities. A mate of a vertex $v$ of a path in $\PP_i$
is on a path of the lower level or the same level with probability $i/K$, thus
\begin{equation}
p_{i}(\o)=1-\frac{i}{K}\;\mbox{.}\label{eq-circ}
\end{equation}
Since the mate is on a path of any fixed level with the same probability equal to $1/K$, we also obtain
\begin{eqnarray}
p_{i}(\+) & = & \sum_{j=1}^{i-1}\frac{1}{K}q^{\o}_{j}(\oplus) \label{eq-plus} \\
p_{i}(\x) & = & \sum_{j=1}^{i-1}\frac{1}{K}q^{\o}_{j}(\otimes) \label{eq-times} \\
p_{i}(\.) & = & \frac{1}{K}+\sum_{j=1}^{i-1}\frac{1}{K}q^{\o}_{j}(\odot) \label{eq-cdot}
\end{eqnarray}
Note that if $i=1$, then $p_{i}(\+)=0$, $p_{i}(\x)=0$ and $p_{i}(\.)=1/K$.

Once we have determined the probabilities of input labels,
we can compute the probabilities of output labels (assuming
that all input labels on a single path are independent of each other). 
If $\RR$ is a set of rules and $\sigma\rightarrow\tau$ is one of the rules,
we use $|\sigma|$ for the length of $\sigma$, $|\tau|$ for the length of $\tau$ (which are the same
in each individual rule) and $\sigma\tau(x,y)$ for the number of symbols $y$ in $\tau$ on positions
of symbols $x$ in $\sigma$ (with $x=\?$ or $y=\?$ being the wild-card),
e.g., if $\sigma=\o\bx\.\+$ and $\tau=\odot\otimes\odot$,
then $|\sigma|=|\tau|=3$, $\sigma\tau(\?,\odot)=2$ and $\sigma\tau(\o,\odot)=1$.
If $\sigma$ is a sequence of input symbols, then $p_i(\sigma)$ denotes the probality that
the sequence $\sigma$ appears on $|\sigma|$ consecutive vertices as their input symbols,
i.e.,
$$p_i(\sigma)=\prod_{j=1}^{|\sigma|}p_i(\sigma_j)\;\mbox{.}$$
Finally, $\norm{\RR}$ is the maximum length of a rule in $\RR$ and
$\overline{\RR}$ is the expansion of $\RR$, the set of rules obtained from $\RR$
by replacing each rule containing the wild-card symbol(s) $\?$
with rules for all possible choices of the value of $\?$,
e.g., a rule containing two symbols $\?$ is replaced with $16$ rules.

The probabilities of output labels of vertices are then given by the following formulas:
\begin{eqnarray}
q_{i}(\oplus) & = & \sum_{\sigma\rightarrow\tau\in\overline{\RR}}\frac{p_i(\sigma)}{Q_i}\sigma\tau(\?,\oplus)
                   \label{eq-oplus} \\
q_{i}(\otimes) & = & \sum_{\sigma\rightarrow\tau\in\overline{\RR}}\frac{p_i(\sigma)}{Q_i}\sigma\tau(\?,\otimes)
                   \label{eq-otimes} \\
q_{i}(\odot) & = & \sum_{\sigma\rightarrow\tau\in\overline{\RR}}\frac{p_i(\sigma)}{Q_i}\sigma\tau(\?,\odot)
                   \label{eq-odot}
\end{eqnarray}
where
$$ Q_i = \sum_{\sigma\rightarrow\tau\in\overline{\RR}} |\sigma|p_i(\sigma)\;\mbox{.}$$
Let us derive (\ref{eq-oplus})--(\ref{eq-odot}) formally.
Consider a path of length $L$ on level $i$.
Let $E$ be the expected number of rules of $\RR$ matched and
$E_{\sigma}$ the expected number of times the rule $\sigma\to\tau$ is matched.
In the analysis that follows,
it is convenient to think of the considered path as a sufficiently long part of an infinite path
which is consistent with our arguments presented later.
When we start matching, the probability that the rule $\sigma\to\tau$ is matched
is equal to $p_i(\sigma)$. Hence,
the expected number of times $E_{\sigma}$ the rule $\sigma\to\tau$ is matched is $p_i(\sigma)E$.
By the linearity of expectations, we also have
that
$$\sum_{\sigma\to\tau\in\overline{\RR}}E_{\sigma}|\sigma|=L\;\mbox{.}$$
By the definition of $Q_i$, we obtain that $E=L/Q_i$ and $E_{\sigma}=Lp_i(\sigma)/Q_i$.
Summing over all rules $\sigma\to\tau\in\overline{\RR}$, we obtain the expected numbers of output labels
of the vertices of a path, e.g., the expected number of the output label $\oplus$ is
$$\sum_{\sigma\rightarrow\tau\in\overline{\RR}}\frac{Lp_i(\sigma)}{Q_i}\sigma\tau(\?,\oplus)\;\mbox{.}$$
The obtained quantities after dividing $L$ represent
the corresponding probabilities as given in (\ref{eq-oplus})--(\ref{eq-odot}).

Similarly, the probabilities conditioned on the appearance of the input symbol $\o$ are given by:
\begin{eqnarray}
q^{\o}_{i}(\oplus) & = & \sum_{\sigma\rightarrow\tau\in\overline{\RR}}\frac{p_i(\sigma)}{Q^{\o}_i}\sigma\tau(\o,\oplus)
                   \label{eq-oplus+} \\
q^{\o}_{i}(\otimes) & = & \sum_{\sigma\rightarrow\tau\in\overline{\RR}}\frac{p_i(\sigma)}{Q^{\o}_i}\sigma\tau(\o,\otimes)
                   \label{eq-otimes+} \\
q^{\o}_{i}(\odot) & = & \sum_{\sigma\rightarrow\tau\in\overline{\RR}}\frac{p_i(\sigma)}{Q^{\o}_i}\sigma\tau(\o,\odot)
                   \label{eq-odot+}
\end{eqnarray}
where
$$ Q^{\o}_i = \sum_{\sigma\rightarrow\tau\in\overline{\RR}} \sigma\tau(\o,\?)p_i(\sigma)\;\mbox{.}$$
The just defined quantities for the set of rules given in Figure~\ref{fig-example} and $K=5$ are given in Figure~\ref{fig-exprob}.

\begin{figure}
\begin{center}
\begin{tabular}{|c|cccc|ccc|ccc|}
\hline
$i$ & $p_{i}(\+)$ & $p_{i}(\x)$ & $p_{i}(\.)$ & $p_{i}(\o)$ & $q_{i}(\oplus)$ & $q_{i}(\otimes)$ & $q_{i}(\odot)$ & $q^{\o}_{i}(\oplus)$ & $q^{\o}_{i}(\otimes)$ & $q^{\o}_{i}(\odot)$ \\
\hline
1 & .0000 & .0000 & .2000 & .8000  &  .0000 & .3333 & .6667 &  .0000 & .3333 & .6667 \\
2 & .0000 & .0667 & .3333 & .6000  &  .0142 & .3160 & .6698 &  .0236 & .3302 & .6462 \\
3 & .0047 & .1327 & .4626 & .4000  &  .0198 & .3009 & .6793 &  .0496 & .3222 & .6282 \\
4 & .0146 & .1972 & .5882 & .2000  &  .0155 & .2889 & .6956 &  .0773 & .3089 & .6138 \\
5 & .0301 & .2590 & .7110 & .0000  &  .0000 & .2812 & .7187 &  - & - & - \\
\hline
\end{tabular}
\end{center}
\caption{The probabilities given by (\ref{eq-circ})--(\ref{eq-odot+}) for the set of rules from Figure~\ref{fig-example} and $K=5$.}
\label{fig-exprob}
\end{figure}

In the ideal case, we understand each of the paths in $\PP_i$ as an infinite path
where the probability that any particular vertex is matched to the first symbol
in one of the rules of $\overline{\RR}$ is the same.
This probability, denoted as $r_{i,1}$, is equal to 
$$\frac{\sum\limits_{\sigma\rightarrow\tau\in\overline{\RR}}p_i(\sigma)}{\sum\limits_{\sigma\rightarrow\tau\in\overline{\RR}}|\sigma|p_i(\sigma)}=\frac{1}{Q_i}\;\mbox{.}$$
This formula can be derived by considering a path of sufficiently large length $L$ and
dividing the expected number of rules matched (which we computed earlier) by $L$.
Similarly, we define $r_{i,j}$ to be the probability that a vertex is matched to the first symbol
of a rule in $\overline{\RR}$ conditioned that none of the $j-1$ previous vertices has been matched to the first symbol
of a rule. Clearly, $r_{i,j}=0$ for $j>\norm{\overline{\RR}}=\norm{\RR}$.
As an example, observe that $r_{1,1}=1/3$, $r_{1,2}=1/2$ and $r_{1,3}=1$ 
for the set of rules from Figure~\ref{fig-example}.

We now describe the actual procedure used to label the vertices of paths and
give all details to cope with technical difficulties which we omit in Subsection~\ref{sub-overview}.
The paths are labelled from the first level, i.e., we start with paths in $\PP_1$, continue
with those in $\PP_2$, etc. In each set $\PP_i$, the paths are labelled in arbitrary order.

Consider a path $P$ from $\PP_i$.
The input label of a vertex $v$ is determined as described in Subsection~\ref{sub-overview}.
Next choose randomly a number $\ell_0$ between $1$ and $\norm{\RR}$ whose value is equal to $\ell$
with probability $r_{i,\ell}\prod_{j=1}^{\ell-1}(1-r_{i,j})$.
Add auxiliary $\norm{\RR}-\ell_0$ vertices at the beginning of $P$ and label each of them with $x$ with probability $p_{i}(x)$ and
add auxiliary $\norm{\RR}$ vertices at the end of $P$ and label each of them with $x$ with probability $p_{i}(x)$.
Finally, apply the labeling procedure as described in Subsection~\ref{sub-overview} (starting
with the $\ell_0$ auxiliary vertices) and when finished, discard the $\ell_0+\norm{\RR}$ auxiliary vertices.
In this way, all vertices of $P$ are assigned output labels.

Let us analyze the actual probabilities that a vertex in $\PP_i$ has a certain input and output labels.
We claim that a vertex of a path $P$ of $\PP_i$ is assigned an input label $x\in\{\+,\x,\.,\o\}$
with probability $p_i(x)$ and an output label $x\in\{\oplus,\otimes,\odot\}$ with probability $q_i(x)$.
Moreover, the probability of a vertex $v$ getting a certain input label depends only on the labels and the levels of 
vertices that are $i$-close to $v$ (and thus the input labels of vertices on $P$ are mutually independent random variables).
Because of the addition of $\norm{\RR}-\ell_0$ auxiliary vertices at the beginning of $P$,
where $\ell_0$ was chosen as described earlier, each vertex of $P$ has the same probability
of being the first, second, etc., in a rule of $\RR$ applied to $P$.
Hence, the probability that its output label is $y\in\{\oplus,\otimes,\odot\}$ is equal to $q_i(y)$.

We now analyze probabilities of input and output labels.
Let $v$ be a vertex of a path $P$ in $\PP_i$.
Since the mate $v'$ of $v$ is on a higher level with probability $1-i/K$,
the input label of $v$ is $\o$ with probability $1-i/K$. With probability $1/K$,
the mate $v'$ is on the same level and the input label of $v$ is $\.$ in this case.
With probability $(i-1)/K$, the mate $v'$ is on a lower level and the input label of $v$
is determined by the output label of $v'$ in this case; the probability that the output label of $v'$
is $y\in\{\oplus,\otimes,\odot\}$ is $q^{\o}_j(y)$ where $j$ is the level of $v'$.
As the output label of $v'$ depends only
on the levels and labels of vertices $(i-1)$-close to other vertices on the path of $v'$,
the output label of $v'$ depends only on the vertices $i$-close to $v$ (and the level of its mate).
Hence, the input label of $v$ depends only on the labels and the levels of vertices $i$-close to $v$.
In particular, input labels of all the vertices of $P$ are mutually independent.
Since the probability of $v$ being the first, second, etc. in a particular rule of $\RR$ applied to it is the same
because of padding with $\ell_0$ auxiliary vertices, the probability that
the output label of $v$ is $y\in\{\oplus,\otimes,\odot\}$ is equal to $q_i(y)$ and
the probability that the output label is $y$ conditioned by its input label being $\o$ is $q^{\o}_i(y)$.

Consider a labeling of the vertices of $G$ constructed in the just described way.
The dominating set $D$ for a graph $G$ is formed by the following vertices:
\begin{itemize}
\item vertices not covered by the paths in $\PP_1\cup\cdots\cup\PP_k$,
\item vertices with input label $\+$,
\item vertices with output label $\otimes$, and
\item the first and the last vertex of each path in $\PP_1\cup\cdots\cup\PP_k$.
\end{itemize}
Let us verify that $D$ is a dominating set: a vertex $v$ not covered by the paths
in $\PP_1\cup\cdots\cup\PP_k$ is in $D$. Vertices on paths $\PP_1\cup\cdots\cup\PP_k$
with input label different from $\+$ and output label different from $\otimes$
are dominated either by their mates or by their neighbors on paths with auxiliary
vertices (assuming the set of rules is correct). However, since the auxiliary
vertices were discarded, the first and the last vertex of each path may not be dominated
in this way---this has been repaired by adding them to $D$ (regardless their
labels).

It remains to estimate the expected size of $D$.
The expected size of $D$ is equal to the sum of
the expected number of vertices with output label $\otimes$ or $\oplus$ (note that
each vertex with input label $\+$ has a mate with output label $\oplus$),
the number of vertices not covered by paths (which is at most $O(n/g)$) and
twice the number of paths because of the inclusion of their first and last vertices to $D$ (this number
is also $O(n/g)$). The probability that a vertex $v$ has output label $\oplus$
is $\frac{q_{1}(\oplus)+\cdots+q_{K}(\oplus)}{K}$ as the path containing $v$ is included to any of the $K$ levels
with the same probability. Similarly, its output label is $\otimes$ with probability
$\frac{q_{1}(\otimes)+\cdots+q_{K}(\otimes)}{K}$.

We summarize the results presented in this section in the following lemma:

\begin{lemma}
\label{lm-main}
Let $\RR$ be a correct set of rules, $K\ge 2$ an integer and $q_{i}(\oplus)$ and $q_{i}(\otimes)$
quantities determined using (\ref{eq-circ})--(\ref{eq-odot+}). If the procedure described in this subsection
given by the set of rules $\RR$ and the integer $K$
is applied to any cubic $n$-vertex graph with girth at least $g$,
then it produces a dominating set with expected size equal to
$$\frac{\sum_{i=1}^K(q_{i}(\oplus)+q_{i}(\otimes))}{K}\;n+O\left(\frac{n}{g}\right)\;\mbox{.}$$
\end{lemma}

Note that already the simple set of rules given in Figure~\ref{fig-example} for $K=5$
yields by Lemma~\ref{lm-main} that
the domination number of a cubic graph with girth at least $g$ is at most $\constfive{}n+O(n/g)$,
an improvement of the bound of Rautenbach and Reed from~\cite{bib-rautenbach+}.

\subsection{Finale}
\label{sub-finale}

We apply Lemma~\ref{lm-main} for a suitable correct set $\RR$ of rules.
This set of $79$ rules can be found in Figure~\ref{fig-rules}; we have
generated this set of rules by hand and checked both ourselves and
by computer that $\RR$ is correct.
For $K=10\,000$, we have computed the values given
in (\ref{eq-circ})--(\ref{eq-odot+}), see Figure~\ref{fig-values}, and
found out that
$$\frac{\sum_{i=1}^K(q_{i}(\oplus)+q_{i}(\otimes))}{K}=\const\;\mbox{.}$$
Lemma~\ref{lm-main} now yields Theorem~\ref{thm-main}.

\begin{figure}
$$\begin{array}{lclclcl}
\bx\+\bx\+\bx\+\bx\? & \rightarrow & \bx\odot\bx\odot\bx\odot\bx\odot
& & \bx\o\bx\.\bx\.\bx\.\bx\.\bx\x & \rightarrow & \bx\otimes\bx\odot\bx\odot\bx\otimes\bx\odot\bx\odot
\\
\bx\+\bx\+\bx\x & \rightarrow & \bx\odot\bx\odot\bx\odot
& & \bx\o\bx\.\bx\.\bx\.\bx\.\bx\. & \rightarrow & \bx\odot\bx\otimes\bx\odot\bx\odot\bx\otimes\bx\odot
\\
\bx\+\bx\+\bx\. & \rightarrow & \bx\odot\bx\odot\bx\odot
& & \bx\o\bx\.\bx\.\bx\.\bx\.\bx\o & \rightarrow & \bx\odot\bx\otimes\bx\odot\bx\odot\bx\otimes\bx\odot
\\
\bx\+\bx\+\bx\o & \rightarrow & \bx\odot\bx\odot\bx\odot
& & \bx\o\bx\.\bx\.\bx\.\bx\o\bx\+\bx\? & \rightarrow & \bx\oplus\bx\odot\bx\otimes\bx\odot\bx\odot\bx\odot\bx\odot
\\
\bx\+\bx\x & \rightarrow & \bx\odot\bx\odot
& & \bx\o\bx\.\bx\.\bx\.\bx\o\bx\x & \rightarrow & \bx\otimes\bx\odot\bx\odot\bx\otimes\bx\odot\bx\odot
\\
\bx\+\bx\. & \rightarrow & \bx\odot\bx\odot
& & \bx\o\bx\.\bx\.\bx\.\bx\o\bx\. & \rightarrow & \bx\odot\bx\otimes\bx\odot\bx\odot\bx\otimes\bx\odot
\\
\bx\+\bx\o & \rightarrow & \bx\odot\bx\odot
& & \bx\o\bx\.\bx\.\bx\.\bx\o\bx\o & \rightarrow & \bx\odot\bx\otimes\bx\odot\bx\odot\bx\otimes\bx\odot
\\
\bx\x & \rightarrow & \bx\odot
& & \bx\o\bx\.\bx\.\bx\o\bx\+\bx\? & \rightarrow & \bx\odot\bx\otimes\bx\odot\bx\odot\bx\odot\bx\odot
\\
\bx\.\bx\x\bx\+\bx\? & \rightarrow & \bx\odot\bx\otimes\bx\odot\bx\odot
& & \bx\o\bx\.\bx\.\bx\o\bx\x & \rightarrow & \bx\oplus\bx\odot\bx\otimes\bx\odot\bx\odot
\\
\bx\.\bx\x\bx\x & \rightarrow & \bx\odot\bx\otimes\bx\odot
& & \bx\o\bx\.\bx\.\bx\o\bx\.\bx\+\bx\? & \rightarrow & \bx\oplus\bx\odot\bx\otimes\bx\odot\bx\odot\bx\odot\bx\odot
\\
\bx\.\bx\x\bx\. & \rightarrow & \bx\odot\bx\otimes\bx\odot
& & \bx\o\bx\.\bx\.\bx\o\bx\.\bx\x & \rightarrow & \bx\otimes\bx\odot\bx\odot\bx\otimes\bx\odot\bx\odot
\\
\bx\.\bx\x\bx\o & \rightarrow & \bx\odot\bx\otimes\bx\odot
& & \bx\o\bx\.\bx\.\bx\o\bx\.\bx\. & \rightarrow & \bx\odot\bx\otimes\bx\odot\bx\odot\bx\otimes\bx\odot
\\
\bx\.\bx\.\bx\+\bx\? & \rightarrow & \bx\odot\bx\otimes\bx\odot\bx\odot
& & \bx\o\bx\.\bx\.\bx\o\bx\.\bx\o & \rightarrow & \bx\odot\bx\otimes\bx\odot\bx\odot\bx\otimes\bx\odot
\\
\bx\.\bx\.\bx\x & \rightarrow & \bx\odot\bx\otimes\bx\odot
& & \bx\o\bx\.\bx\.\bx\o\bx\o\bx\+\bx\? & \rightarrow & \bx\oplus\bx\odot\bx\otimes\bx\odot\bx\odot\bx\odot\bx\odot
\\
\bx\.\bx\.\bx\. & \rightarrow & \bx\odot\bx\otimes\bx\odot
& & \bx\o\bx\.\bx\.\bx\o\bx\o\bx\x & \rightarrow & \bx\otimes\bx\odot\bx\odot\bx\otimes\bx\odot\bx\odot
\\
\bx\.\bx\.\bx\o & \rightarrow & \bx\odot\bx\otimes\bx\odot
& & \bx\o\bx\.\bx\.\bx\o\bx\o\bx\. & \rightarrow & \bx\odot\bx\otimes\bx\odot\bx\odot\bx\otimes\bx\odot
\\
\bx\.\bx\o\bx\+\bx\? & \rightarrow & \bx\odot\bx\otimes\bx\odot\bx\odot
& & \bx\o\bx\.\bx\.\bx\o\bx\o\bx\o & \rightarrow & \bx\odot\bx\otimes\bx\odot\bx\odot\bx\otimes\bx\odot
\\
\bx\.\bx\o\bx\x & \rightarrow & \bx\odot\bx\otimes\bx\odot
& & \bx\o\bx\.\bx\o\bx\+\bx\? & \rightarrow & \bx\otimes\bx\odot\bx\odot\bx\odot\bx\odot
\\
\bx\.\bx\o\bx\. & \rightarrow & \bx\odot\bx\otimes\bx\odot
& & \bx\o\bx\.\bx\o\bx\x & \rightarrow & \bx\odot\bx\otimes\bx\odot\bx\odot
\\
\bx\.\bx\o\bx\o & \rightarrow & \bx\odot\bx\otimes\bx\odot
& & \bx\o\bx\.\bx\o\bx\.\bx\+\bx\? & \rightarrow & \bx\odot\bx\otimes\bx\odot\bx\odot\bx\odot\bx\odot
\\
\bx\.\bx\+\bx\+\bx\? & \rightarrow & \bx\odot\bx\odot\bx\odot\bx\odot
& & \bx\o\bx\.\bx\o\bx\.\bx\x & \rightarrow & \bx\oplus\bx\odot\bx\otimes\bx\odot\bx\odot
\\
\bx\.\bx\+\bx\x & \rightarrow & \bx\odot\bx\odot\bx\odot
& & \bx\o\bx\.\bx\o\bx\.\bx\.\bx\+\bx\? & \rightarrow & \bx\oplus\bx\odot\bx\otimes\bx\odot\bx\odot\bx\odot\bx\odot
\\
\bx\.\bx\+\bx\. & \rightarrow & \bx\odot\bx\odot\bx\odot
& & \bx\o\bx\.\bx\o\bx\.\bx\.\bx\x & \rightarrow & \bx\otimes\bx\odot\bx\odot\bx\otimes\bx\odot\bx\odot
\\
\bx\.\bx\+\bx\o & \rightarrow & \bx\odot\bx\odot\bx\odot
& & \bx\o\bx\.\bx\o\bx\.\bx\.\bx\. & \rightarrow & \bx\odot\bx\otimes\bx\odot\bx\odot\bx\otimes\bx\odot
\\
\bx\o\bx\x & \rightarrow & \bx\oplus\bx\odot
& & \bx\o\bx\.\bx\o\bx\.\bx\.\bx\o & \rightarrow & \bx\odot\bx\otimes\bx\odot\bx\odot\bx\otimes\bx\odot
\\
\bx\o\bx\+\bx\+\bx\? & \rightarrow & \bx\odot\bx\odot\bx\odot\bx\odot
& & \bx\o\bx\.\bx\o\bx\.\bx\o\bx\+\bx\? & \rightarrow & \bx\oplus\bx\odot\bx\otimes\bx\odot\bx\odot\bx\odot\bx\odot
\\
\bx\o\bx\+\bx\x & \rightarrow & \bx\odot\bx\odot\bx\odot
& & \bx\o\bx\.\bx\o\bx\.\bx\o\bx\x & \rightarrow & \bx\otimes\bx\odot\bx\odot\bx\otimes\bx\odot\bx\odot
\\
\bx\o\bx\+\bx\. & \rightarrow & \bx\odot\bx\odot\bx\odot
& & \bx\o\bx\.\bx\o\bx\.\bx\o\bx\. & \rightarrow & \bx\odot\bx\otimes\bx\odot\bx\odot\bx\otimes\bx\odot
\\
\bx\o\bx\+\bx\o & \rightarrow & \bx\odot\bx\odot\bx\odot
& & \bx\o\bx\.\bx\o\bx\.\bx\o\bx\o & \rightarrow & \bx\odot\bx\otimes\bx\odot\bx\odot\bx\otimes\bx\odot
\\
\bx\o\bx\o\bx\+\bx\? & \rightarrow & \bx\oplus\bx\odot\bx\odot\bx\odot
& & \bx\o\bx\.\bx\o\bx\o\bx\+\bx\? & \rightarrow & \bx\odot\bx\otimes\bx\odot\bx\odot\bx\odot\bx\odot
\\
\bx\o\bx\o\bx\x & \rightarrow & \bx\oplus\bx\oplus\bx\odot
& & \bx\o\bx\.\bx\o\bx\o\bx\x & \rightarrow & \bx\oplus\bx\odot\bx\otimes\bx\odot\bx\odot
\\
\bx\o\bx\o\bx\. & \rightarrow & \bx\odot\bx\otimes\bx\odot
& & \bx\o\bx\.\bx\o\bx\o\bx\.\bx\+\bx\? & \rightarrow & \bx\oplus\bx\odot\bx\otimes\bx\odot\bx\odot\bx\odot\bx\odot
\\
\bx\o\bx\o\bx\o & \rightarrow & \bx\odot\bx\otimes\bx\odot
& & \bx\o\bx\.\bx\o\bx\o\bx\.\bx\x & \rightarrow & \bx\otimes\bx\odot\bx\odot\bx\otimes\bx\odot\bx\odot
\\
\bx\o\bx\.\bx\+\bx\? & \rightarrow & \bx\oplus\bx\odot\bx\odot\bx\odot
& & \bx\o\bx\.\bx\o\bx\o\bx\.\bx\. & \rightarrow & \bx\odot\bx\otimes\bx\odot\bx\odot\bx\otimes\bx\odot
\\
\bx\o\bx\.\bx\x & \rightarrow & \bx\otimes\bx\odot\bx\odot
& & \bx\o\bx\.\bx\o\bx\o\bx\.\bx\o & \rightarrow & \bx\odot\bx\otimes\bx\odot\bx\odot\bx\otimes\bx\odot
\\
\bx\o\bx\.\bx\.\bx\+\bx\? & \rightarrow & \bx\otimes\bx\odot\bx\odot\bx\odot\bx\odot
& & \bx\o\bx\.\bx\o\bx\o\bx\o\bx\+\bx\? & \rightarrow & \bx\oplus\bx\odot\bx\otimes\bx\odot\bx\odot\bx\odot\bx\odot
\\
\bx\o\bx\.\bx\.\bx\x & \rightarrow & \bx\odot\bx\otimes\bx\odot\bx\odot
& & \bx\o\bx\.\bx\o\bx\o\bx\o\bx\x & \rightarrow & \bx\otimes\bx\odot\bx\odot\bx\otimes\bx\odot\bx\odot
\\
\bx\o\bx\.\bx\.\bx\.\bx\+\bx\? & \rightarrow & \bx\odot\bx\otimes\bx\odot\bx\odot\bx\odot\bx\odot
& & \bx\o\bx\.\bx\o\bx\o\bx\o\bx\. & \rightarrow & \bx\odot\bx\otimes\bx\odot\bx\odot\bx\otimes\bx\odot
\\
\bx\o\bx\.\bx\.\bx\.\bx\x & \rightarrow & \bx\oplus\bx\odot\bx\otimes\bx\odot\bx\odot
& & \bx\o\bx\.\bx\o\bx\o\bx\o\bx\o & \rightarrow & \bx\odot\bx\otimes\bx\odot\bx\odot\bx\otimes\bx\odot
\\
\bx\o\bx\.\bx\.\bx\.\bx\.\bx\+\bx\? & \rightarrow & \bx\oplus\bx\odot\bx\otimes\bx\odot\bx\odot\bx\odot\bx\odot
\end{array}$$
\caption{The set of $79$ rules used in the proof of Theorem~\ref{thm-main}.}
\label{fig-rules}
\end{figure}

\begin{figure}
\begin{center}
\begin{tabular}{|c|cccc|ccc|ccc|}
\hline
$i$ & $p_{i}(\+)$ & $p_{i}(\x)$ & $p_{i}(\.)$ & $p_{i}(\o)$ & $q_{i}(\oplus)$ & $q_{i}(\otimes)$ & $q_{i}(\odot)$ & $q^{\o}_{i}(\oplus)$ & $q^{\o}_{i}(\otimes)$ & $q^{\o}_{i}(\odot)$ \\
\hline
      1 & .0000 & .0000 & .0001 & .9999  &  .0000 & .3333 & .6667  &  .0000 & .3333 & .6667 \\
      2 & .0000 & .0033 & .0002 & .9998  &  .0000 & .3333 & .6667  &  .0000 & .3333 & .6666 \\
 \vdots & \vdots& \vdots& \vdots& \vdots &  \vdots& \vdots& \vdots &  \vdots& \vdots& \vdots \\
 2\,500 & .0094 & .0808 & .1598 & .7500  &  .0542 & .2853 & .6605  &  .0723 & .3165 & .6112 \\
 \vdots & \vdots& \vdots& \vdots& \vdots &  \vdots& \vdots& \vdots &  \vdots& \vdots& \vdots \\
 5\,000 & .0352 & .1593 & .3055 & .5000  &  .0661 & .2480 & .6859  &  .1322 & .3131 & .5547 \\
 \vdots & \vdots& \vdots& \vdots& \vdots &  \vdots& \vdots& \vdots &  \vdots& \vdots& \vdots \\
 7\,500 & .0744 & .2382 & .4373 & .2500  &  .0447 & .2223 & .7329  &  .1790 & .3207 & .5003 \\
 \vdots & \vdots& \vdots& \vdots& \vdots &  \vdots& \vdots& \vdots &  \vdots& \vdots& \vdots \\
 9\,999 & .1222 & .3211 & .5566 & .0001  &  .0000 & .2106 & .7894  &  .1967 & .3465 & .4568 \\
10\,000 & .1222 & .3211 & .5566 & .0000  &  .0000 & .2106 & .7894  &  - & - & - \\
\hline
\end{tabular}
\end{center}
\caption{The numerical values of probabilities given by (\ref{eq-circ})--(\ref{eq-odot+}) for the set of rules from Figure~\ref{fig-rules} and $K=10\,000$.}
\label{fig-values}
\end{figure}

We have also constructed (with computer assistance)
a correct set of $3607$ rules such that Lemma~\ref{lm-main}
applied for $K=1\,000\,000$ yields that every cubic graph
has domination number at most $\constsup\;n+O(n/g)$ but
we neither present nor claim this bound here.

\section*{Acknowledgement}

This research was conducted while the authors visited University of Ljubljana.
The authors would like to thank Riste {\v S}krekovski who was hosting them
for an extra-ordinary working conditions and for numerous discussions on the subject.

\pagebreak

\end{document}